

\baselineskip=14pt
\parskip=10pt

\font\eightrm=cmr8 
\font\eighttt=cmtt8
\magnification=\magstephalf

\def\1{{\overline{1}}}
\def\2{{\overline{2}}}
\parindent=0pt
\overfullrule=0in
\def\Tilde{\char126\relax}
\def\frac#1#2{{#1 \over #2}}
\bf
\centerline
{
The $C$-FINITE ANSATZ
}
\rm
\bigskip
\centerline{ {\it
Doron 
ZEILBERGER}\footnote{$^1$}
{\eightrm  \raggedright
Department of Mathematics, Rutgers University (New Brunswick),
Hill Center-Busch Campus, 110 Frelinghuysen Rd., Piscataway,
NJ 08854-8019, USA.
{\eighttt zeilberg  at math dot rutgers dot edu} ,
\hfill \break
{\eighttt http://www.math.rutgers.edu/\~{}zeilberg/} .
Written: July 16, 2011.
Accompanied by Maple package {\eighttt Cfinite}
downloadable from 
{\eighttt http://www.math.rutgers.edu/\~{}zeilberg/tokhniot/Cfinite} .
Supported in part by the NSF.
}
}
\bigskip
\centerline
{\it Dedicated with friendship and admiration to Mourad Ismail and Dennis Stanton}
\bigskip
{\bf Apology}

While this article is dedicated to both Mourad Ismail and Dennis Stanton, it does not directly reference
any of their works. The main reason is that I only talk about the most {\it trivial} kind of recurrences:
{\bf linear} and {\bf constant coefficients}. But try searching {\tt Ismail AND Recurrences}
or {\tt Stanton AND Recurrences} in the database {\tt MathSciNet} (or {\tt Google Scholar})
and you would see that both Mourad and Dennis are great {\it gurus}
in recurrences, so the subject matter of this paper is not entirely inappropriate as a tribute to them. 
The present work is also largely {\it experimental}, and Dennis Stanton is a great pioneer in computer experimentations!

{\bf PROLOGUE}

Before starting the paper itself, let me very briefly mention what I talked 
about at the wonderful conference

{\it q-Series 2011: An International Conference on q-Series, Partitions and Special Functions},

``honouring {\bf Mourad Ismail} and 
{\bf Dennis Stanton} for their valuable contributions to Number theory and Special Functions throughout their careers'',
that took place on March 14-16, 2011 at Georgia Southern University, and perfectly organized
by {\bf Drew Sills}. My {\bf plenary} talk was entitled

``{\it Some Golden Oldies of Mourad Ismail and Dennis Stanton}''

and the abstract was short and sweet:

{\it ``This is the time to recall some of the beautiful mathematics that I learned from Dennis and Mourad''}.

I talked, {\it inter alia}, about Dennis Stanton's {\bf amazing} article

``{\it A short proof of a generating function for Jacobi polynomials}'',

that appeared in
{\it  Proc. Amer. Math. Soc.}  {\bf 80}  (1980), 398-400, where he used an idea that he attributed to
his advisor, Dick Askey, but that Askey modestly claims goes back to Hermite, to give a very
{\it soft} and elegant proof of a very hairy formula of Bailey.

Another thing that I talked about was Mourad Ismail's {\bf gorgeous} article 

``{\it A simple proof of Ramanujan's ${}_1\psi_1$ sum}'',

that contained the {\it proof-from-the-book} of {\bf Ramanujan's} lovely formula
$$
\sum_{n=-\infty}^{\infty} \frac{(a)_n}{(b)_n}t^n \, = \,
\frac{(b/a)_{\infty}(at)_{\infty}(q/at)_{\infty}(q)_{\infty}}
{(q/a)_{\infty}(b/at)_{\infty}(b)_{\infty}(t)_{\infty}} \quad ,
$$
(where, as usual $(A)_n:=(1-A) \cdots (1-q^{n-1}A)$, if $n \geq 0$,
$(A)_n:=1/((1-q^{-1}A) \cdots (1-q^nA))$, if $n <0$, 
and $(A)_{\infty}:=(1-A)(1-qA) (1-q^2A)\cdots$ ) .

Mourad's article appeared in the {\it proceedings} (of the AMS) three years earlier (PAMS {\bf  63} (1977), 185-186.)
This proof was already {\it immortalized} and {\it canonized} (and {\it shrunk} to half a page!)
in Appendix C of George Andrews' classic monograph, {\bf CBMS \#66}, that appeared in 1986 and was based on ten 
beautiful lectures delivered at an NSF-CBMS conference that took place in 1985, at Arizona State University,
and organized by {\bf Mourad Ismail} and Ed Ihring.

In my talk, I mentioned that while Mourad's proof  certainly qualifies to be included  in ``God's book'',
since, like Mourad and George, God is, by definition, an {\it infinitarian}, it does not qualify
to be included in {\it my} book. It said (in George's rendition, my {\bf emphasis})

``Regarding  the left-hand side  as an {\bf analytic} function of $b$ for $|b|<1$ \dots '' .

Then one plugs-in $b=q^N$ ($N=0,1,2,\dots$), getting the trivial {\bf q-binomial theorem}, and one
sees that the left side minus the right side vanishes for ``infinitely'' many values of $b$ and
then uses the ``fact'' that an ``analytic'' function inside 
$|b|<1$ that vanishes
on a ``convergent'' ``infinite'' sequence ``must'' be identically zero.

Of course, to {\bf finitists} like myself, this proof is entirely non-rigorous, since it uses
{\it fictional} things like so-called analytic functions, and uses heavy guns from a sophisticated
(and flawed!) ``infinite'' theory.

But don't despair! It is very easy to translate Mourad's flawed proof and make it entirely legit,
and in the process make it even nicer. Replace the phrase

{\it ``analytic function of $b$ defined in $|b|<1$''}

by

{\it ``bilateral formal power series in $t$ whose coefficients are rational functions of $b$''} ,

and note that the difference of the left and right sides is a bilateral formal power series (in $t$) whose
coefficients are rational functions of $b$ (and of course also of $a$ and $q$ but that's irrelevant).
A rational function of $b$, whose degree of the numerator is, say, $m$, is identically zero if it vanishes at 
$m+1$ distinct values of $b$, so the ``infinitely'' many points $b=q^N$ ($N=0,1,2 \dots$) is
more than enough.

{\bf FEATURE PRESENTATION}

{\bf The Maple package} {\tt Cfinite} 

This article is accompanied by  the Maple package {\tt Cfinite} downloadable directly from

{\eighttt http://www.math.rutgers.edu/\~{}zeilberg/tokhniot/Cfinite},

or from the ``front'' 

{\eighttt http://www.math.rutgers.edu/\~{}zeilberg/mamarim/mamarimhtml/cfinite.html},

where one can find links to
fiftenn sample input and output files, some of which are mentioned
throughout this paper.

{\bf $C$-finite Sequences} 

Recall that a $C$-finite sequence $\{a(n)\}, n=0,1,\dots$ is a sequence that satisfies a
{\bf linear-recurrence equation with  constant coefficients}. It is known (but not as well-known as it should be!)
and easy to see (e.g. [Z2],[KP]) that the set of $C$-finite sequences is an {\it algebra}.
Even though a $C$-finite sequence is an ``infinite'' sequence, it is in fact, {\bf like everything else in mathematics}
(and elsewhere!) a {\bf finite} object. 
An order-$L$ $C$-finite sequence $a(n)$
is completely specified by the coefficients $c_1,c_2, \dots , c_L$  of the recurrence
$$
a(n)=c_1 a(n-1)+ c_2 a(n-2) + \dots + c_L a(n-L) \quad,
$$
and the {\bf initial conditions}
$$
a(0)=d_1 \quad , \quad \dots \quad , \quad a(L-1)=d_L \quad.
$$
So a $C$-finite sequence can be {\bf coded} in terms of the $2L$ ``bits'' of information
$$
[[d_1, \dots , d_L],[c_1, \dots, c_L] ] \quad .
$$
For example,  the Fibonacci sequence is written:
$$
[[0,1],[1,1]] .
$$
Since this ansatz (see [Z2]) is fully decidable, it is
possible to decide equality, and  evaluate {\it ab initio}, wide classes of sums, and
things are easier than the {\it holonomic ansatz}[Z1]. The wonderful new book
by Manuel Kauers and Peter Paule[KP] 
also presents a convincing case. See [HX][GW][Kau][KZ]
for very interesting and efficient algorithms.

{\bf Rational Generating Functions}

Equivalently, a $C$-finite sequence is a sequence $\{a(n)\}$ whose {\bf ordinary generating function},
$\sum_{n=0}^{\infty} a(n)z^n$,
is a {\bf rational function} of $z$. These come up a lot in combinatorics and elsewhere
(e.g. formal languages). See the {\it old testament}[St2], chapter 4,  and the {\it new testament}[KP], chapter 4.

[To go from a $C$-finite  representation to a rational function, use {\tt CtoR(C,z);} in {\tt Cfinite}.
To go the other way, do {\tt RtoC(f,z);} .]

{\bf Etymology}

I coined the term {\it $C$-finite sequence} in [Z1], 
as a {\it hybrid} analog of  Richard Stanely's [St1] names ``$D$-finite function''
and ``$P$-recursive sequence''. If I had to do it over I would call them ``{\it $C$-recursive sequences}'',
but it is too late now
since the term $C$-finite  already
made it into the wonderful {\it undergraduate} textbook [KP], and it is also in the title of the important paper [GW].

{\bf Zeilberger-style proofs: You (Often) CAN generalize from FINITELY Many Cases }

The conventional wisdom of mathematics (at least for the last 2500 years), preached to us by our teachers and that, in turn, we preach to
our students, is that {\it you can't generalize from finitely many cases}. While this is certainly true 
{\it sometimes}, it is not {\it always} true. Many times you {\bf can} generalize from  {\bf finitely} many cases, just
like natural scientists.

Michael Hirschhorn kindly called this style of proof
{\it ``in the spirit of Zeilberger''}, see his beautiful proof[H] of an amazing identity of Ramanujan, that gives {\it infinitely}
many ``almost'' counterexamples to Fermat's Last Theorem for $n=3$, namely infinitely many
triples $\{ (a,b,c) \} $ such that $a^3+b^3+c^3=\pm 1$.

{\bf But not Everyone Knows About this Style of Proof}

Everybody knows that {\it numerical} identities like $2+7=3 \times 3$ are {\it routinely provable}, using {\it standard algorithms}.
But even  people like Neil Sloane (and Jeff Shallit) and James Sellers are not fully
aware that identities amongst $C$-finite sequences are equally routine.
If they were, James Sellers' article [Se] would not have been accepted for publication in {\it Journal of Integer Sequences},
or written in the first place, see the parody [Sr].   
There are hundreds (possibly thousands) of articles like this in the literature, sometimes giving
``elegant'' proofs of such trivial results. While it is always nice to have elegant proofs,
honesty requires that the authors state clearly, in the abstract, that the result that they are
elegantly proving is routinely provable. See my opinion [Z4].

{\bf The First Reason for This Article: Educating}

Since even Neil Sloane, Jeff Shallit, James Sellers, Art Benjamin, and many other people are not
(fully) aware of the {\it triviality} of the $C$-finite ansatz (or more politely,
there being an {\it algorithmic proof theory} for it), and in spite of the articles and book cited above,
I thought that it is a good idea to make it better known.

{\bf The Second Reason for This Article: Implementation}

While Curtis Greene and Herb Wilf [GW] and  Manuel Kauers[Kau] (see also [KZ]) already have Mathematica implementations
of many operations on $C$-finite sequences,
and possibly also Maple  ones, I thought that it is a good idea to design
a $C$-finite calculator, that also enables one to {\it discover} new identities.
The novelty, in that part, is the
{\it approach}, pure guessing! (that is justified {\it a posteriori}).

{\bf The Third Reason for This Article: Factorization}

The truly novel part (I believe) is in addressing the problem of 
{\bf factorization}. See below.

{\bf The $C$-finite Calculator}

In order to {\it decipher} a $C$-finite sequence where the first few terms are given,
all you need is use linear algebra to ``guess'' (using the {\bf ansatz}) the $c$'s (you already know the $d$'s).
See Procedure {\tt GuessRec} in {\tt Cfinite}.
If you have two $C$-finite sequences $C_1$ and $C_2$ of order $L_1$ and $L_2$, you don't
need any fancy footwork to figure out the sequence $C_1+C_2$ (of order $L_1+L_2$)
and the sequence $C_1C_2$ (of order $L_1L_2$). All you need is to crank out
$L_1+L_2+4$ (the $+4$ is for safety reasons) and $L_1L_2+4$ terms, respectively, and let the
computer guess the $C$-finite description, completely by {\bf guessing},
using {\it undetermined coefficients} that is implemented by Procedure {\tt GuessRec} in {\tt Cfinite}.

So it is (very!) easy to {\bf multiply} $C$-finite sequences, in other words, from the
$C$-finite descriptions of $C_1$ and $C_2$ get the $C$-finite description of $C_1C_2$
(by $C_1C_2$ we mean the sequence whose $n$-th term is $C_1(n)C_2(n)$,
in terms of their generating functions it is called the  {\it Hadamard product}).

{\bf Procedures } {\tt Khibur} and {\tt Kefel}  of {\tt Cfinite}

In order to add two C-finite sequences $C1$ and $C2$, simply type

{\tt Khibur(C1,C2);} \quad .

In order to multiply type:

{\tt Kefel(C1,C2);} \quad ,

and for a verbose version, presenting fully detailed {\it proofs in the spirit of Zeilberger}, type \hfill\break
{\tt KefelV(C1,C2);} .

{\bf An Example of Using Procedure} {\tt Kefel}

Let $\{U_n(x)\}_{n=0}^{\infty}$ be the sequence of Chebyshev polynomials of the second kind, i.e. the sequence of polynomials in $x$ defined by
$$
\sum_{n=0}^{\infty} U_n(x)t^n \, = \, \frac{1}{1-2xt+t^2} \quad .
$$

Typing

{\tt CtoR(Kefel(Ux(a),Ux(b)),t);} would give you, in a few seconds, the following result[Sh]

{\bf Lou Shapiro's Product-Of-Two-Chebyshev-Polynomials Identity}
$$
\sum_{n=0}^{\infty} U_n(a) U_n(b) \, t^n \, = \, {\frac {1-{t}^{2}}{1-4\,abt- \left( -4\,{a}^{2}+2-4\,{b}^{2} \right) {t}^{2}-4\,ab{t}^{3}+{t}^{4}}} \quad.
$$

Typing

{\tt  CtoR(Kefel(Kefel(Ux(a),Ux(b)),Ux(c)),t);}

yields, in a few more seconds, the following much deeper result

{\bf Shalosh B. Ekhad's Product-Of-Three-Chebyshev-Polynomials Identity}
$$
\sum_{n=0}^{\infty} \, U_n(a) U_n(b) U_n(c) t^n \, = \frac{N(t)}{D(t)} \quad ,
$$
where the polynomials $N(t)$ are $D(t)$ are as follows.

$$
N(t)=
1+ \left( -4\,{a}^{2}-4\,{b}^{2}-4\,{c}^{2}+3 \right) {t}^{2}+16\,abc{t
}^{3}+ \left( -4\,{a}^{2}-4\,{b}^{2}-4\,{c}^{2}+3 \right) {t}^{4}+{t}^{
6} \quad ,
$$
$$
D(t)=
{t}^{8}-8\,abc{t}^{7}+ \left( 16\,{a}^{2}{b}^{2}+16\,{a}^{2}{c}^{2}-8
\,{a}^{2}+16\,{b}^{2}{c}^{2}-8\,{b}^{2}-8\,{c}^{2}+4 \right) {t}^{6}+
 \left( -32\,{a}^{3}bc+40\,abc-32\,a{b}^{3}c-32\,ab{c}^{3} \right) {t}^
{5}+ 
$$
$$
\left( 16\,{a}^{4}+64\,{a}^{2}{b}^{2}{c}^{2}-16\,{a}^{2}+16\,{b}^{
4}-16\,{b}^{2}+6+16\,{c}^{4}-16\,{c}^{2} \right) {t}^{4}+ 
\left( -32\,{a}^{3}bc+40\,abc-32\,a{b}^{3}c-32\,ab{c}^{3} \right) {t}^{3}
$$
$$
+ \left( 16
\,{a}^{2}{b}^{2}+16\,{a}^{2}{c}^{2}-8\,{a}^{2}+16\,{b}^{2}{c}^{2}-8\,{b
}^{2}-8\,{c}^{2}+4 \right) {t}^{2}-8\,abct +1 \quad .
$$

{\bf The (Computationally) Hard Problem of Factoring  C-finite Sequences}

Alas, going backwards 
(just like in integer factorization, that makes our ATM cards hopefully secure)
is not so easy! If you are given a $C$-finite sequence of order $L$, say, and $L$ is  composite,
$L=L_1 L_2$, (with $L_1,L_2 >1$) you would like to know whether there exist $C$-finite sequence
$C_1$ and $C_2$ such that $C=C_1C_2$, and if they do, find them.
One way, that works for small $L$, is to
do {\it symbolic} multiplication of {\it generic} $C$-finite sequences, and then try to solve,
by matching coefficients, the resulting {\bf non-linear} system of algebraic equations, using
the Buchberger algorithm. [This is implemented  in procedure {\tt Factorize} of the Maple package {\tt Cfinite} .]
But for larger orders this is hopeless! 

Procedure {\tt FactorizeI1} does the same by {\it brute force}, but 
only handles {\it integer} sequences.
While it can't go very far, it {\it discovered}, {\it ab initio}, in less than half a second, the three factorizations in [Se].
See

{\tt http://www.math.rutgers.edu/\Tilde zeilberg/tokhniot/oCfinite4} \quad .

For a more verbose version see

{\tt http://www.math.rutgers.edu/\Tilde zeilberg/tokhniot/oCfinite5} \quad ,

and for a completely spelled-out proofs, {\it in the spirit of Zeilberger}, see

{\tt http://www.math.rutgers.edu/\Tilde zeilberg/tokhniot/oCfinite5a} \quad .

{\bf Why is this problem interesting?}

Two great landmarks of Statistical Physics are the Onsager[O] solution of the 
two-dimensional Ising model and
the Kasteleyn[Kas]- Temperley-Fisher[TF] solutions of the dimer problem.
They use lots of human (ad-hoc) ingenuity to first 
get an {\bf explicit} answer for a {\it finite strip} of 
{\it arbitrary} (symbolic)  width.
They then take the so-called {\it thermodynamic limit}. 
It turns out that in either case, the $m$-wide strip sequence is a $C$-finite sequence of order $2^m$.
Surprisingly, in both cases they happen to be products of $m$ $C$-finite sequences of order $2$
(different, but closely related).

Since nowadays computers can automatically, {\bf completely rigorously},
figure out the $C$-finite description for each strip of width $m$,
for {\it specific}, numeric $m$,
([EZ],[Z3])
(in practice easily for $m \leq 10$), knowing how to ``factorize'' them explicitly, would
lead one to conjecture both solutions, with fully rigorous proofs for $m \leq 10$, and
with larger computers, beyond. Since physicists are not as hung-up as mathematicians
about {\it rigorous} proofs, that would have been a great breakthrough, even {\it without}
the {\it human} proofs for general $m$. Besides, the explicit ``conjecture'' discovered by
the computer might suggest and inspire (to obtuse mathematicians) a formal proof.

{\bf The ``Cheating Algorithm''}

Since it is so  hard to factorize explicitly, it is still nice to know, as fast as possible,
whether or not the inputted $C$-finite sequence is factorizable. If it is not, it would be stupid
to waste efforts in trying to factorize it. If it is, then it is worthwhile applying for
time on a bigger computer.

Recall (Binet) that ``generically'' every $C$-finite sequence
$[[d_1, \dots , d_L],[c_1, \dots, c_L] ]$
of order $L$ can be written as a linear combination
$$
C(n)=\sum_{i=1}^{L} C_i \alpha_i^n \quad,
$$
where the $C_i$'s depend on the initial conditions, and the $\alpha_i$'s are the roots of the
{\it characteristic equation}
$$
z^L-\sum_{i=1}^{L} c_iz^{L-i}=0 \quad .
$$

So if $C:=C_1C_2$ and the roots of $C_1$ and $C_2$ are
$\alpha_1, \dots, \alpha_{L_1}$, and
$\beta_1, \dots, \beta_{L_2}$ respectively, then the roots of $C$,
let's call them $\gamma_1, \dots, \gamma_{L_1L_2}$, 
consist
of the {\bf Cartesian product}
$$
\{ \alpha_i  \beta_j \, | \, 1 \leq i \leq L_1 \quad , \quad 1 \leq j \leq L_2 \} \quad .
$$
with $L:=L_1 L_2$ elements.

If this is indeed the case, then the set of $L^2$ ratios 
$$
\{ \frac{\gamma_i}{\gamma_j}\, | \, 1 \leq i , j \leq L \} \quad ,
$$
would have a certain
{\it profile of repetitions} that the computer can easily figure out for
{\it arbitrary} symbols $\alpha_1, \dots, \alpha_{L_1}$ and
$\beta_1, \dots, \beta_{L_2}$.  

{\bf Procedure} {\tt ProdIndictor} of the Maple package {\tt Cfinite}

Procedure {\tt ProdIndicator(m,n)} yields the {\it profile of repetitions}
indicative of the characteristic roots of  a C-finite
sequence that happens to be the product of a C-finite sequence of order $m$ and a C-finite sequence
of order $n$.

For example, {\tt ProdIndicator(2,2)} yields:
$$
 [1, 1, 1, 1, 2, 2, 2, 2, 4] \quad .
$$
To understand what is going on, let's work it out {\it by hand}.
$$
\gamma_1=\alpha_1 \beta_1 \quad , \quad \gamma_2= \alpha_1 \beta_2  \quad , \quad \gamma_3= \alpha_2 \beta_1 
\quad , \quad \gamma_4= \alpha_2 \beta_2  \quad .
$$
So
$$
\frac{\gamma_1}{\gamma_1}=\frac{\alpha_1 \beta_1}{\alpha_1 \beta_1}  =1 
\quad, \quad
\frac{\gamma_1}{\gamma_2}=\frac{\alpha_1 \beta_1}{\alpha_1 \beta_2}=  \frac{\beta_1}{ \beta_2}
\quad, \quad
\frac{\gamma_1}{\gamma_3}=\frac{\alpha_1 \beta_1}{\alpha_2 \beta_1}=  \frac{\alpha_1}{\alpha_2}
\quad, \quad
\frac{\gamma_1}{\gamma_4}=\frac{\alpha_1 \beta_1}{\alpha_2 \beta_2}  \quad ,
$$
$$
\frac{\gamma_2}{\gamma_1}=\frac{\alpha_1 \beta_2}{\alpha_1 \beta_1}  =\frac{\beta_2}{ \beta_1}  
\quad, \quad
\frac{\gamma_2}{\gamma_2}=\frac{\alpha_1 \beta_2}{\alpha_1 \beta_2}  =1
\quad, \quad
\frac{\gamma_2}{\gamma_3}=\frac{\alpha_1 \beta_2}{\alpha_2 \beta_1}
\quad, \quad
\frac{\gamma_2}{\gamma_4}=\frac{\alpha_1 \beta_2}{\alpha_2 \beta_2}=  \frac{\alpha_1}{\alpha_2} \quad ,
$$
$$
\frac{\gamma_3}{\gamma_1}=\frac{\alpha_2 \beta_1}{\alpha_1 \beta_1} =\frac{\alpha_2}{\alpha_1} 
\quad, \quad
\frac{\gamma_3}{\gamma_2}=\frac{\alpha_2 \beta_1}{\alpha_1 \beta_2}
\quad, \quad
\frac{\gamma_3}{\gamma_3}=\frac{\alpha_2 \beta_1}{\alpha_2 \beta_1}  =1
\quad, \quad
\frac{\gamma_3}{\gamma_4}=\frac{\alpha_2 \beta_1}{\alpha_2 \beta_2}  = \frac{\beta_1}{\beta_2}   \quad ,
$$
$$
\frac{\gamma_4}{\gamma_1}=\frac{\alpha_2 \beta_2}{\alpha_1 \beta_1}  
\quad, \quad
\frac{\gamma_4}{\gamma_2}=\frac{\alpha_2 \beta_2}{\alpha_1 \beta_2}=\frac{\alpha_2 }{\alpha_1}
\quad, \quad
\frac{\gamma_4}{\gamma_3}=\frac{\alpha_2 \beta_2}{\alpha_2 \beta_1}=\frac{\beta_2}{\beta_1}
\quad, \quad
\frac{\gamma_4}{\gamma_4}=\frac{\alpha_2 \beta_2}{\alpha_2 \beta_2} =1 \quad .
$$
We see that the multi-set of all $16$ ratios has:
{\bf four} occurrences of $1$,
{\bf two} occurrences each of 
$\frac{\beta_1}{\beta_2}$,$\frac{\beta_2}{\beta_1}$,
$\frac{\alpha_1}{\alpha_2}$,$\frac{\alpha_2}{\alpha_1}$, and 
four {\bf singletons},
namely $\frac{\alpha_2 \beta_2}{\alpha_1 \beta_1}$,  $\frac{\alpha_2 \beta_1}{\alpha_1 \beta_2}$
and their reciprocals.

Now for the proposed $C$-finite sequence of order $L$, find
(in floating point!, but with {\tt Digits:=100;}) approximations to the roots of its characteristic equation,
then form these $L^2$ ratios, and  group them into classes with the ``same'' value
(up to the agreed-on approximation). If you get the same {\it pattern of repetition}, then you have
proved (empirically) that the given $C$-finite sequence $C$, of order $L=L_1L_2$, is indeed the product of
$C$-finite sequences of orders $L_1$ and $L_2$.
Procedure {\tt IsProd}  in {\tt Cfinite} implements this algorithm. See the source code for more details.

If $C$ has order $L=L_1 L_2 \cdots L_r$, and you want to find out whether $C$ is a product of
$r$ $C$-finite sequences of orders 
$L_1, \dots, L_r$ you do the analogous thing.
Procedure {\tt IsProdG}  in {\tt Cfinite} implements this 
more general scenario.

{\bf Output}

Using the output from [EZ] we confirmed that the straight enumeration dimer 
problems for strips of width
$\leq 10$ are indeed products of $C$-finite sequences of order $2$. See \hfill\break
{\tt http://www.math.rutgers.edu/\Tilde zeilberg/tokhniot/oCfinite2} \quad  .

Using the output from [Z3] we confirmed that the weighted enumeration dimer problem for strips of width
$\leq 10$ indeed are products of $C$-finite sequences of order $2$ for many random numerical assignments
of the weights.

As for the actual factorization, we were, on our modest computer, only able to find them (from scratch, without peeking at the answer) for
$m \leq 7$, but a more clever implementation, and a larger computer, no doubt would
be able to {\bf conjecture} {\bf ab initio} (without any human ad-hocery!) the
exact solution of the dimer problem derived and proved in [Kas] and [TF].
Ditto for Onsager's [O] (human) {\it tour-de-force}.

Sample output for some of the other procedures (e.g {\tt BT} for the Binomial Transform
and {\tt GuessNLR} for finding  non-linear (polynomial) recurrences of lower-order than the
(linear) order of a given $C$-finite sequence) can be obtained from the front of the present
article:

  {\tt http://www.math.rutgers.edu/\Tilde zeilberg/mamarim/mamarimhtml/cfinite.html} .

{\bf Encore: 1142 beautiful and deep Greene-Wilf-style Fibonacci identities in less than 4400 seconds}

See the computer-generated webbook 

{\tt http://www.math.rutgers.edu/\Tilde zeilberg/tokhniot/oCfinite13} \quad  .

For simpler identities see also 

{\tt http://www.math.rutgers.edu/\Tilde zeilberg/tokhniot/oCfinite11}  and 

{\tt http://www.math.rutgers.edu/\Tilde zeilberg/tokhniot/oCfinite12} .

These fascinating new identities can keep bijective combinatorialists busy for the next one hundred years.
Each of these identities cries out for an insightful, elegant, combinatorial proof!

\vfill\eject

{\bf References}

[EZ] Shalosh B. Ekhad and Doron Zeilberger,
{\it Automatic Generation of Generating Functions for Enumerating Matchings}, 
Personal Journal of Shalosh B. Ekhad
and Doron Zeilberger, April 29, 2011,
{\tt http://www.math.rutgers.edu/\Tilde zeilberg/mamarim/mamarimhtml/shidukhim.html} \quad .

[GW] Curtis Greene and Herbert S. Wilf,
{\it Closed form summation of $C$-finite sequences}, Trans. Amer. Math. Soc. {\bf 359}(2007), 1161-1189,
{\tt http://arxiv.org/abs/math/0405574} \quad .

[H] M.D. Hirschhorn, {\it A proof in the spirit of Zeilberger of an amazing identity of Ramanujan}, 
Math. Mag. {\bf 4 }(1996), 267-269.

[HX] Qing-Hu Hou and  Guoce Xin
{\it Constant term evaluation for summation of C-finite sequences}, DMTCS proc. AN, 2010, 761-772.

[Kas] P. W. Kasteleyn, {\it The statistics of dimers on a lattice: I. The number of
dimer arrangements in a quadratic lattice}, Physica {\bf 27} (1961), 1209-1225.

[Kau] Manuel Kauers,
{\it SumCracker: A package for manipulating symbolic sums and related objects},
Journal of Symbolic Computation {\bf 41} (2006),  1039-1057,
\hfill\break
{\tt http://www.risc.uni-linz.ac.at/people/mkauers/publications/kauers06h.pdf} \quad .

[KP] Manuel Kauers  and Peter Paule, {\it ``The Concrete Tetrahedron''}, Springer, 2011.

[KZ] Manuel Kauers and Burkhard Zimmermann,
{\it  Computing the algebraic relations of C-finite sequences and multisequences},
Journal of Symbolic Computation {\bf 43} (2008), 787-803.

[O] Lars Onsager, {\it Crystal statistics, I. A two-dimensional model with an order-disorder transition},
Phys. Rev. {\bf 65}(1944), 117-149.

[Se] James A. Sellers, {\it Domino Tilings and Products of Fibonacci and Pell Numbers},
Journal of Integer Sequences, {\bf 5}(2002), 02.1.2, 
\hfill\break
{\tt http://www.cs.uwaterloo.ca/journals/JIS/VOL5/Sellers/sellers4.pdf} \quad .

[Sh] Louis Shapiro, {\it A combinatorial proof of a Chebyshev polynomial identity}, Discrete Math. {\bf 34}(1981), 203-206 .

[Sr] Semaj Srelles, {\it The Sum of Two and Seven and The Product of Three by Three },
Personal J. of Shalosh B. Ekhad and Doron Zeilberger, April 12, 2011,
\hfill\break
{\tt http://www.math.rutgers.edu/\Tilde zeilberg/mamarim/mamarimhtml/parodia.html} \quad .

[St1] Richard Stanley, {\it Differentiably finite power series}, European Journal of Combinatorics {\bf 1}(1980), 175-188,
{\tt http://www-math.mit.edu/~rstan/pubs/pubfiles/45.pdf } \quad .

[St2] Richard Stanley, { \it ``Enumerative combinatorics''}, volume {\bf 1}, Wadsworth and Brooks/Cole, Pacific Grove, CA, 1986,
second printing, Cambridge University Press, Cambridge, 1996.

[TF] H. Temperley and M. Fisher, {\it Dimer Problems in Statistical
Mechanics-an exact result}, Philos. Mag. {\bf 6} (1961), 1061-1063.

[Z1] Doron Zeilberger, {\it  A Holonomic Systems Approach To Special Functions},
J. Computational and Applied Math {\bf 32}(1990), 321-368,
\hfill\break
{\tt http://www.math.rutgers.edu/\Tilde zeilberg/mamarim/mamarimPDF/holonomic.pdf } \quad .

[Z2] Doron Zeilberger,
{\it An Enquiry Concerning Human (and Computer!) [Mathematical] Understanding}
in: C.S. Calude ,ed., ``Randomness and Complexity, from Leibniz to Chaitin'', World Scientific, Singapore, Oct. 2007,
\hfill\break
{\tt http://www.math.rutgers.edu/\Tilde zeilberg/mamarim/mamarimhtml/enquiry.html} \quad .

[Z3] Doron Zeilberger, {\it Automatic CounTilings}, Personal Journal of Shalosh B. Ekhad
and Doron Zeilberger, Jan. 20, 2006; also published in Rejecta Mathematica
{\bf 1} (2009), 10 -17,
\hfill\break
{\tt http://www.math.rutgers.edu/\Tilde zeilberg/mamarim/mamarimhtml/tilings.html}  \quad .

[Z4] Doron Zeilberger,
{\it Opinion 89: Mental Math Whiz [And Very Good Mathematician] Art Benjamin Should be 
Aware that not only his Night Job, but also some parts of his ``Day Job'' should be clearly labeled ``For Entertainment Only'' },
\hfill\break
{\tt http://www.math.rutgers.edu/\Tilde zeilberg/Opinion89.html} \quad .

\end